\documentclass[10pt]{amsart}
\usepackage{amsfonts}
\usepackage{color}
\usepackage[square,compress,comma, numbers]{natbib}
\usepackage[colorlinks=true, citecolor=blue, linkcolor=blue]{hyperref}
\allowdisplaybreaks[4]
\usepackage{amssymb}
\usepackage{amsmath}
\usepackage{color}

\definecolor{c20}{rgb}{0.,0.7,0.}
\definecolor{c30}{rgb}{0.,0.,1.}
\definecolor{c40}{rgb}{1,0.1,0.7}
\definecolor{c50}{rgb}{1,0,0}
\definecolor{c60}{rgb}{1,0.9,0.1}

\allowdisplaybreaks[4]

\newcommand{\abs}[1]{\left\lvert #1 \right\rvert}

\newcommand{\E}[1]{\mathbb{E}\left\{#1\right\}}

\newcommand{\pk}[1]{\mathbb{P} \left \{#1 \right \} }
\newcommand{\R}{\mathbb{R}}

\newcommand{\BQN}{\begin{eqnarray}}
\newcommand{\EQN}{\end{eqnarray}}
\newcommand{\BQNY}{\begin{eqnarray*}}
\newcommand{\EQNY}{\end{eqnarray*}}

\newcommand{\BS}{\begin{sat}}
\newcommand{\ES}{\end{sat}}
\newcommand{\BT}{\begin{theo}}
\newcommand{\ET}{\end{theo}}
\newcommand{\BL}{\begin{lem}}
\newcommand{\EL}{\end{lem}}
\newcommand{\BK}{\begin{korr}}
\newcommand{\EK}{\end{korr}}

\newcommand{\BD}{\begin{de}}
\newcommand{\ED}{\end{de}}
\newcommand{\BIT}{\begin{itemize}}
\newcommand{\EIT}{\end{itemize}}
\newcommand{\BDI}{\begin{description}}
\newcommand{\EDI}{\end{description}}

\newcommand{\BRM}{\begin{remarks}}
\newcommand{\ERM}{\end{remarks}}

\newcommand{\BEL}{\begin{lem}}
\newcommand{\EEL}{\end{lem}}

\newtheorem{theo}{Theorem}[section]
\newtheorem{sat}[theo]{Proposition}
\newtheorem{de}[theo]{Definition}
\newtheorem{lem}[theo]{Lemma}

\newtheorem{korr}[theo]{Corollary}

\newtheorem{remarks}[theo]{Remarks}

\newcommand{\prooftheo}[1]{ \textsc{\bf Proof of Theorem} \ref{#1}:}

\newcommand{\COM}[1]{}

\newcommand{\QED}{\hfill $\Box$}

\def\rw{\rightarrow}

\def\IF{\infty}

\def\LT{\left}
\def\RT{\right}

\def\rw{\rightarrow}

\def\Var{\text{Var}}

\def\Bu+#1{\mathcal{B}^{\varepsilon+}_{u}(#1)}

\begin{document}

\title{Extremes of multifractional Brownian motion}

\author{Long Bai}
\address{Long Bai,
Department of Mathematical Sciences, Xi'an Jiaotong-Liverpool University,
Suzhou 215123
China, and Department of Actuarial Science, University of Lausanne, UNIL-Dorigny, 1015 Lausanne, Switzerland
}
\email{long.bai@xjtlu.edu.cn}
\bigskip

 \maketitle

{\bf Abstract:}
Let $B_{H}(t), t\geq [0,T], T\in(0,\infty)$ be the standard Multifractional Brownian Motion(mBm), in this contribution we are concerned with the exact asymptotics of
\begin{eqnarray*}
 \mathbb{P}\left\{\sup_{t\in[0,T]}B_{H}(t)>u\right\}
\end{eqnarray*}
as $u\rightarrow\infty$. Mainly depended on the structures of  $H(t)$, the results under several important cases are investigated.


\section{Introduction}

Multifractional Brownian motion generalize the fractional fractional Brownian motion(fBm)of exponent $H$ to the case where $H$ is no longer a constant but a function of the time index of the process.  It was introduced to overcome certain limitations of the classical fBm. Contrarily to fBm, the almost sure H\"{o}lder exponent of mBm is allowed to vary along the trajectory, a useful feature when one needs to model processes whose regularity evolves in time, such as internet traffic or images.
The tail behaviors of the supremum of the fBm over certain time interval have been investigated in the literature, such as \cite{ANTGFM2002, ITAFRM2011, ECCGP1999, SGRPFBM2013}. But there are rare results about the extreme of mBm. In this paper we try to extend the exact tail asymptotics to mBm case, i.e.  we consider for $0\leq T_1<T_2<\IF$
\BQN\label{pkBH}
\pk{\sup_{t\in[T_1,T_2]}B_{H}(t)>u},\quad u\rw\IF.
\EQN

Through this paper, $H(t): [0,\IF)\rw[\underline{H},\overline{H}]\subset(0,1)$, the index function of mBm, is a positive H\"{o}lder function of exponent $\lambda>0$, i.e.
\BQN\label{Ht1}
\abs{H(t)-H(s)}\leq \mathbb{C} \abs{t-s}^\lambda, t,s\in[0,T],
\EQN
where $\mathbb{C}$ is a positive constant and satisfies
\BQN\label{Ht2}
0<\forall_{t\geq 0}H(t)\leq \min(1,\lambda).
\EQN
First we give the definition of standard mBm $B_{H}(t), t\geq 0$ as in \cite{mBm1995}.
\BD
Let $Y_H(t), t\geq 0$ be a mBm and $H(t)$ satisfy \eqref{Ht1} and \eqref{Ht2}. Then there exists a unique continuous positive function $C(H(t)), t\geq 0$ such that the process $B_{H}(t), t\geq 0$ defined by $B_{H}(t)=\frac{Y(t)}{C(H(t))}$ is continuous and verifies the following property
\BQNY
\Var\LT(\frac{B_{H}(t+h)-B_{H}(t)}{h^{H(t)}}\RT)\rw 1,\quad  h\rw 0.
\EQNY
The process $B_{H}(t), t\geq 0$ is called Standard Multifractional Brownian Motion. \ED
Then by \cite{EGRP1997} and \cite{mBmA2000}, we have the following representation of the standard mBm.
\BD
{\bf(Standard mBm, Harmonizable Representation)}\\
Let $H(t)$  satisfy \eqref{Ht1} and \eqref{Ht2}. For $t\geq 0$ the following function is called standard mBm:
\BQNY
B_{H}(t)=\frac{1}{C(H(t))}\int_{\R}\frac{e^{it\xi}-1}{\abs{\xi}^{H(t)+\frac{1}{2}}}d B(\xi)
\EQNY
where $B$ denotes standard Brownian motion and
\BQNY
C(x)=\LT(\frac{\pi}{x\Gamma(2x)\sin(\pi x)}\RT)^{1/2}.
\EQNY
\ED
\section{Main results}
We introduce a lemma which is an important property of the mBm and some notation before we state our main results.
\BEL
By \cite{mBmA2000}, the explicit expressions for the autocovariance  of $B_{H}(t),t\geq 0$ is
\BQNY
&&Cov(B_{H}(t),B_{H}(s))=\E{B_{H}(t)B_{H}(s)}=D(H(s),H(t))\LT(s^{H(s)+H(t)}+t^{H(s)+H(t)}-|t-s|^{H(s)+H(t)}\RT),
\EQNY
where $$D(x,y)=\frac{\sqrt{\Gamma(2x+1)\Gamma(2y+1)\sin(\pi x)\sin(\pi y)}}{2\Gamma(x+y+1)\sin\LT(\frac{\pi(x+y)}{2}\RT)}.$$

Further by $D(H(t),H(t))=\frac{1}{2}$
\BQNY
\sigma^2_{H}(t):=\E{B^2_{H}(t)}=t^{2H(t)}, \ t\geq 0.
\EQNY
By \cite{atlocally}, the correlation function of ${B}_{H}(t)$ satisfies for any $t> 0$
\BQN\label{rr}
r_{H}(t,t+h)=1-\frac{1}{2}t^{-2H(t)}|h|^{2H(t)}+o(|h|^{2H(t)}),\ \ h\rw 0.
\EQN
\EEL

Next we need to introduce some notation, starting with the well-known Pickands constant $\mathcal{H}_\alpha$ defined by
$$\mathcal{H}_{\alpha}=\lim_{S\rightarrow\IF}\frac{1}{S}\mathcal{H}_{\alpha}[0,S],\quad \text{with }
\mathcal{H}_{\alpha}[0,S]=\E{\sup_{t\in[0,S]}e^{\sqrt{2}B_\alpha(t)-|t|^\alpha}}\in(0,\IF),
$$
where $S>0$ is a constant and $B_\alpha(t),t\in\mathbb{R}$ is a standard fractional Brownian motion (fBm)  with Hurst index $\alpha/2\in(0,1].$
Further, define for $a>0$
\BQNY
\mathcal{P}_{\alpha}^{a}=\lim_{S\rw\IF}\mathcal{P}_{\alpha}^{a}[0,S],\quad
\widetilde{\mathcal{P}_{\alpha}^{a}}=\lim_{S\rw\IF}\mathcal{P}_{\alpha}^{a}[-S,S],
\EQNY
with
\BQNY
\mathcal{P}_{\alpha}^{a}[-S,S]=
\E{\sup_{t\in[-S,S]}e^{\sqrt{2}B_\alpha(t)-|t|^\alpha-a\abs{t}^\alpha}}.
\EQNY
These constants defined above play a significant role in the following theorems,
see \cite{Tabis} for various properties of these constants and compare with, e.g.,
\cite{GeneralPit16,PicandsA,Pit72, DE2002,DI2005,DE2014,DiekerY,DEJ14,Pit20, DM, SBK, Htilt,DHL14Ann,DebKo2013}.

Now we return to our principal problem deriving below the exact asymptotic behaviour of \eqref{pkBH}.
\BT\label{Thm1}
 If $\sigma_{H}(t)$, the standard deviation of $B_H(t)$, attains its maximum over $[0,T]$ at  $t^*\in [0,T]$, and $H(t)$ satisfies that
\BQN\label{ASH}
H(t^*+h)=H(t^*)-c\abs{h}^\gamma+o(\abs{h}^\gamma),\quad h\rw 0,\quad c>0, \quad \gamma\in(0,1).
\EQN
Then we have as $u\rw\IF$
\BQNY
\pk{\sup_{t\in[0,T]}B_{H}(t)>u}\sim \Psi\LT(\mu\RT)
\LT\{
\begin{array}{ll}
\LT(1+\mathbb{I}_{\{t^*\in(0,T)\}}\RT)\mathcal{H}_{\alpha}a^{1/\alpha}
b^{-1/\gamma}\Gamma\LT(\frac{1}{\gamma}+1\RT)\mu^{\frac{2}{\alpha}-\frac{2}{\gamma}},& \ \text{if}\ \ \alpha<\gamma,\\
\widehat{\mathcal{P}^{b/a}_{\alpha}},&\  \text{if}\ \ \alpha=\gamma,\\
1,&\ \text{if}\ \ \alpha>\gamma,
\end{array}
\RT.
\EQNY
where
\BQNY
\widehat{\mathcal{P}^{b/a}_{\alpha}}=\LT\{
\begin{array}{ll}
\mathcal{P}^{b/a}_{\alpha}, &\ \text{if}\ t^*\in\{0,T\},\\
\widetilde{\mathcal{P}^{b/a}_{\alpha}}, &\ \text{if}\ t^*\in(0,T),
\end{array}
\RT.
\EQNY
 $\mu=u{t^*}^{-H(t^*)}$, $a=\frac{1}{2}{t^*}^{-2H(t^*)}$, $b=c\ln t^*$, and $\alpha=2H(t^*)$.
\ET
Following are several special cases which can not be included in the last theorem.
\BT\label{Thm2}
i) If
\BQNY
H(t)=\frac{1}{\ln t},\quad t\in[T_1,T_2],\quad e<T_1<T_2<\IF,
\EQNY
we have as $u\rw\IF$
\BQNY
\pk{\sup_{t\in[T_1,T_2]}B_{H}(t)>u}\sim \mathcal{H}_{\alpha}a^{1/\alpha}b^{-1}
\frac{\mu^{2/\alpha}}{\ln \mu}\Psi\LT(\mu\RT),
\EQNY
where $\mu=\frac{u}{e}$, $a=\frac{1}{2}{T_2}^{-2H(T_2)}$, $b=\frac{2}{\alpha^2T_2(\ln T_2)^2}$, and $\alpha=2H(T_2)$.\\
ii) If
$$H(t)=ct^{\gamma}, t\in[T_1,T_2], 0<T_1<T_2<\IF, \ c, \gamma>0,$$
we  have as $u\rw\IF$
\BQNY
\pk{\sup_{t\in[T_1,T_2]}B_{H}(t)>u}\sim \Psi\LT(\mu\RT)
\LT\{
\begin{array}{ll}
Q\mathcal{H}_{\alpha}a^{1/\alpha}b^{-1}\mu^{\frac{2}{\alpha}-2},& \ \text{if}\ \ \alpha<1,\\
Q\mathcal{P}^{b/a}_{\alpha},&\  \text{if}\ \ \alpha=1,\\
Q,&\ \text{if}\ \ \alpha>1,
\end{array}
\RT.
\EQNY
where
\BQNY
Q=\LT\{
\begin{array}{ll}
2, &\ \text{if}\ \sigma_H(T_1)=\sigma_H(T_2),\\
1, &\ \text{other},
\end{array}
\RT.
\quad\quad
\widetilde{T}=\LT\{
\begin{array}{ll}
T_1, &\ \text{if}\ \sigma_H(T_1)>\sigma_H(T_2),\\
T_2, &\ \text{if}\ \sigma_H(T_1)\leq\sigma_H(T_2),
\end{array}
\RT.
\EQNY $\mu=u\widetilde{T}^{-c\widetilde{T}^\gamma}$, $a=\frac{1}{2}{\widetilde{T}}^{-2H(\widetilde{T})}$, $b=c \widetilde{T}^{\gamma-1}\abs{1+\gamma\ln \widetilde{T}}$, and $\alpha=2H(\widetilde{T})$.\\
iii) If $H(t)$ is a differentiable function and decreases over $[T_1,T_2]$ with $ T_2\leq1$ or increases over $[T_1,T_2]$ with $T_2\geq1$,
we have as $u\rw\IF$
\BQNY
\pk{\sup_{t\in[T_1,T_2]}B_{H}(t)>u}\sim \Psi\LT(\mu\RT)
\LT\{
\begin{array}{ll}
\mathcal{H}_{\alpha}a^{1/\alpha}b^{-1}\mu^{\frac{2}{\alpha}-2},& \ \text{if}\ \ \alpha<1,\\
\mathcal{P}^{b/a}_{\alpha},&\  \text{if}\ \ \alpha=1,\\
1,&\ \text{if}\ \ \alpha>1,
\end{array}
\RT.
\EQNY
where $\mu=uT_2^{-H(T_2)}$, $a=\frac{1}{2}{T_2}^{-2H(T_2)}$, $b=\frac{H(T_2)}{T_2}+H'(T_2)\ln T_2$ and $\alpha=2H(T_2)$.
\ET
\section{Proofs}
In this section, we give the proofs of all the results.\\
\prooftheo{Thm1} By \eqref{ASH}, we have
\BQNY
\sigma_H(t^*+h)&=&e^{H(t^*+h)\ln(t^*+h)}\\
&=&e^{\LT(H(t^*)-c\abs{h}^\gamma+o(\abs{h}^\gamma)\RT)\LT(\ln t^*+\ln(1+\frac{h}{t^*})\RT)}\\
&=&e^{H(t^*)\ln t^*}
\LT(1+\frac{H(t^*)}{t^*}h\ln\LT(1+\frac{h}{t^*}\RT)^{\frac{t^*}{h}}-c(\ln t^*)\abs{h}^\gamma+o(\abs{h}^\gamma)\RT)\\
&=&{t^*}^{H(t^*)}
\LT(1-c(\ln t^*)\abs{h}^\gamma+o(\abs{h}^\gamma)\RT), h\rw 0.
\EQNY
 Thus by \eqref{rr}, $B_{H}(t)$ is a centered Gaussian case with the standard deviation function and correlation function given as
\BQNY
\sigma_H(t^*+h)= {t^*}^{H(t^*)}
\LT(1-c\ln t^*\abs{h}^\gamma+o(\abs{h}^\gamma)\RT), \quad h\rw 0,
\EQNY
and
\BQNY
r_{H}(t,s)=1-\frac{1}{2}{t^*}^{-2H(t^*)}|s-t|^{2H(t^*)}+o(|s-t|^{2H(t^*)}),\ \ s,t\rw t^*.
\EQNY
Further, since $\sigma_{H}(t)$ attains its maximum over $[0,T]$ at a unique point  $t^*\in [0,T]$, by \cite{GauTrend16} [Proposition 3.9] the result follows.

\QED

\prooftheo{Thm2} i) We have
\BQNY
\sigma_H(t)=e^{H(t)\ln t}=e,\ t\in[T_1,T_2].
\EQNY
$H(t)$ attains its minimum over $[T_1,T_2]$ at $T_2$ with
\BQNY
H(T_2-h)= \frac{1}{\ln T_2}+\frac{1}{T_2(\ln T_2)^2} h+o(h^2),\quad h\downarrow 0,
\EQNY
which combined with \eqref{rr}, we have that $B_H(t)$ is a $\alpha(t)$-locally stationary Gaussian processes as in  \cite{atlocally}. Thus the results follows from  \cite{atlocally} [Theorem 2.1].\\
ii) We have
$H'(t)=c\gamma t^{\gamma-1}$ and
 \BQNY
 \sigma'_H(t)&=&e^{H(t)\ln t}\LT(H'(t)\ln t+\frac{H(t)}{t}\RT)\\
 &=&t^{ct^\gamma}\LT(c\gamma t^{\gamma-1}\ln t+c t^{\gamma-1}\RT), \quad t\in[T_1,T_2]\subset(0,\IF).
 \EQNY
 Setting $\sigma'_H(t)=0$, we get $\widetilde{t}=e^{-\frac{1}{\gamma}}$ and $\sigma_H(t)$ is decreasing on $(0,\widetilde{t})$ and increasing on $(\widetilde{t},\IF)$.
 Thus $\sigma_H(t)$ attains its maximum over $[T_1,T_2]$ at $T_1$ or $T_2$.\\
Case 1: When $\sigma_H(T_1)>\sigma_H(T_2)$, we have $T_1<e^{-\frac{1}{\gamma}}$ and
$$ H(T_1+h)=c T_1^\gamma+c\gamma T_1^{\gamma-1} h+o(h),\ h\downarrow 0.$$
Further,
\BQNY
\sigma_H(T_1+h)&=&e^{H(T_1+h)\ln(T_1+h)}\\
&=&e^{ \LT(cT_1^\gamma+c\gamma T_1^{\gamma-1} h+o(h) \RT)\LT(\ln T_1+\ln(1+\frac{h}{T_1})\RT)}\\
&=&e^{cT_1^\gamma \ln T_1}
\LT(1+c T_1^{\gamma-1}h\ln\LT(1+\frac{h}{T_1}\RT)^{\frac{T_1}{h}}+c\gamma T_1^{\gamma-1}h\ln T_1 +o(h)\RT)\\
&=&T_1^{cT_1^\gamma}
\LT(1+\LT(c T_1^{\gamma-1}+c\gamma T_1^{\gamma-1}\ln T_1\RT)h +o(h)\RT)\\
&=&T_1^{cT_1^\gamma}
\LT(1-c T_1^{\gamma-1}\abs{1+\gamma\ln T_1}h +o(h)\RT),\ h\downarrow 0,
\EQNY
which combined with \eqref{rr} and \cite{GauTrend16} [Proposition 3.9] derives the result.\\
Case 2: When $\sigma_H(T_1)<\sigma_H(T_2)$, we have $\sigma_H(t)$ attains its maximum over $[T_1,T_2]$ at $T_2$,
 \BQNY
 H(T_2-h)=c T_2^\gamma-c\gamma T_2^{\gamma-1} h+o(h),\ h\downarrow 0,
 \EQNY
and
\BQNY
\sigma_H(T_2-h)&=&e^{H(T_2-h)\ln(T_2-h)}\\
&=&e^{ \LT(cT_2^\gamma-c\gamma T_2^{\gamma-1} h+o(h) \RT)\LT(\ln T_2+\ln(1-\frac{h}{T_2})\RT)}\\
&=&e^{cT_2^\gamma \ln T_2}
\LT(1-c T_2^{\gamma-1}h\ln\LT(1-\frac{h}{T_2}\RT)^{-\frac{T_2}{h}}-c\gamma T_2^{\gamma-1}h\ln T_2 +o(h)\RT)\\
&=&T_2^{cT_2^\gamma}
\LT(1-c T_2^{\gamma-1}\abs{1+\gamma\ln T_2}h +o(h)\RT),\ h\downarrow 0,
\EQNY
which combined with \eqref{rr} and \cite{GauTrend16} [Proposition 3.9] derives the result.\\
Case 3: If $\sigma_H(T_1)=\sigma_H(T_2)$, i.e., $\sigma_H(t)$ attains its maximum over $[T_1,T_2]$ at two point $T_1$ and $T_2$, by \cite{Pit96}[Corollary 8.2], we have for $\delta\in(0,\frac{T_2-T_1}{2})$
\BQNY
\pk{\sup_{t\in[T_1,T_2]}B_{H}(t)>u}\sim\pk{\sup_{t\in[T_1,T_1+\delta]}B_{H}(t)>u}
+\pk{\sup_{t\in[T_2-\delta,T_2]}B_{H}(t)>u}, \ u\rw\IF.
\EQNY
Then for $\pk{\sup_{t\in[T_1,T_1+\delta]}B_{H}(t)>u}$ and $\pk{\sup_{t\in[T_2-\delta,T_2]}B_{H}(t)>u}$, using the similar argument as in Case 1 and Case 2, we get the asymptotics.\\
iii) When $H(t)$ is a decrease function over $[T_1,T_2],$ with $T_2\leq1$, we have
 \BQNY
 \sigma'_H(t)&=&e^{H(t)\ln t}\LT(H'(t)\ln t+\frac{H(t)}{t}\RT)>0, \quad t\in[T_1,T_2]\subset(0,1].
 \EQNY
When $H(t)$ is a increase function over $[T_1,T_2],$ with $T_2\geq1$, we have $\sigma_H(t)<1=\sigma_H(1)$ for $t\in(0,1)$ and
 \BQNY
 \sigma'_H(t)&=&e^{H(t)\ln t}\LT(H'(t)\ln t+\frac{H(t)}{t}\RT)>0, \quad t\in[1,\IF).
 \EQNY
Thus we have in both cases $\sigma_H(t)$ attains its maximum over $[T_1,T_2]$ at $T_2$ and
 \BQNY
 H(T_2-h)=H(T_2)-H'(T_2)h+o(h),\ h\downarrow 0.
 \EQNY
 We have
\BQNY
\sigma_H(T_2-h)&=&e^{H(T_2-h)\ln(T_2-h)}\\
&=&e^{ \LT(H(T_2)-H'(T_2)h+o(h)\RT)\LT(\ln T_2+\ln(1-\frac{h}{T_2})\RT)}\\
&=&e^{H(T_2)\ln T_2}
\LT(1-\frac{H(T_2)h}{T_2}\ln\LT(1-\frac{h}{T_2}\RT)^{-\frac{T_2}{h}}-H'(T_2)h\ln T_2 +o(h)\RT)\\
&=&T_2^{H(T_2)}
\LT(1-\LT(\frac{H(T_2)}{T_2}+H'(T_2)\ln T_2\RT)h +o(h)\RT),\ h\downarrow 0,
\EQNY
which combined with \eqref{rr} and \cite{GauTrend16} [Proposition 3.9] derives the result.

\QED

{\bf Acknowledgement}: Thanks to Professor Enkelejd Hashorva who give many useful suggestions which greatly improve our manuscript. Thanks to  Swiss National Science Foundation grant no.  200021-175752.

\bibliographystyle{ieeetr}
\bibliography{mBm}
\end{document}